\documentclass[11pt,a4paper]{article}
\usepackage{amsmath,amssymb,amsthm,geometry}
\usepackage[T1]{fontenc}
\usepackage[utf8]{inputenc}
\theoremstyle{plain}
\newtheorem{theorem}{Theorem}[section]
\newtheorem{proposition}[theorem]{Proposition}
\newtheorem{corollary}[theorem]{Corollary}
\begin{document}
\title{A Matrix-Valued Quotient and a Compatible Bilinear Product on Euclidean Vector Spaces}
\author{J. Enrique H. Ramirez}
\maketitle
\begin{abstract}
For vectors $a,b\in\mathbb R^n$, $b\ne0$, we study a matrix-valued quotient that maps $b$ to $a$ and a compatible bilinear matrix-valued product. The quotient is $E(a,b)=\|b\|^{-2}[(a\cdot b)I+ab^{\mathsf T}-ba^{\mathsf T}]$; the associated product is $M(a,b)=\|b\|^2E(a,b)$. We establish transport and covariance identities, derive a unique vector inverse from the matrix identity $M(a,x)=I$, and characterize the product by orthogonal equivariance, transport, and either a scalar symmetric-part normalization or a diagonal normalization. A constrained Frobenius-norm principle characterizes the quotient. Its spectral decomposition gives exact criteria for invertibility and shows why matrix reciprocity fails in dimensions greater than two. The construction is placed in the context of complex multiplication and exterior and Clifford algebra, without claiming historical priority for its constituent formulas.
\end{abstract}
\noindent\textbf{Keywords: }matrix-valued quotient; bilinear map; orthogonal equivariance; skew-symmetric operator; isotropic tensor; exterior algebra

\noindent\textbf{MSC 2020: }15A69; 15A24; 15A72

\section{Introduction}
Let $n\ge2$. For $a,b\in\mathbb R^n$ with $b\ne0$, the matrix equation $Xb=a$ has infinitely many solutions. The rank-one operator $ab^{\mathsf T}/\|b\|^2$ is one solution, and adding any operator that annihilates $b$ produces another. Consequently, the transport condition alone does not define a quotient of vectors.

An earlier preprint by Ram'irez and Roca~\cite{ramirezroca2023} introduced matrix-valued vector division and an associated multiplication in dimensions two and three. The present article builds on that initial construction, extending the formulation to arbitrary finite-dimensional Euclidean spaces and focusing on axiomatic uniqueness, a Frobenius-norm characterization, and spectral properties. The previously introduced operations are not presented here as new.

This paper examines a specific selection: an operator whose symmetric part is scalar and whose skew-symmetric part is explicitly $ab^{\mathsf T}-ba^{\mathsf T}$ up to normalization. The resulting matrix-valued quotient is compatible with a bilinear product. We establish the algebraic relations between these operations, determine the vector inverse by solving a matrix equation, and characterize the construction by equivariance and a constrained minimization property. The relation to standard constructions involving complex numbers, bivectors, and Clifford algebra is discussed in Section~\ref{sec:context} \cite{hestenes,clifford}.

Here the word \emph{quotient} refers to a matrix satisfying a prescribed transport identity. It does not imply that the matrix is nonsingular or that scalar division laws extend unchanged to arbitrary dimensions.

\section{Matrix-valued quotient and product}
Write $I=I_n$ and use the Euclidean inner product. For $a,b\in\mathbb R^n$ with $b\ne0$, define
\begin{equation}\label{eq:E}
 E(a,b)=\frac{(a\cdot b)I+ab^{\mathsf T}-ba^{\mathsf T}}{\|b\|^2}.
\end{equation}
For arbitrary $a,b\in\mathbb R^n$, define the bilinear map
\begin{equation}\label{eq:M}
 M(a,b)=(a\cdot b)I+ab^{\mathsf T}-ba^{\mathsf T}.
\end{equation}
The codomain of $M$ is $\mathbb R^{n\times n}$, rather than $\mathbb R^n$; hence $M$ is not an internal binary product on the vector space.

\begin{proposition}\label{prop:basic}
For all admissible arguments,
\begin{equation}\label{eq:basic}
 E(a,b)b=a,\qquad M(a,b)b=\|b\|^2a,
\end{equation}
\begin{equation}
 E(a,a)=I\quad(a\ne0),\qquad M(a,a)=\|a\|^2I,
\end{equation}
and
\begin{equation}\label{eq:transpose}
 M(b,a)=M(a,b)^{\mathsf T},\qquad E(a,b)=\frac{M(a,b)}{\|b\|^2}.
\end{equation}
For every orthogonal $Q\in O(n)$,
\begin{equation}\label{eq:cov}
 M(Qa,Qb)=QM(a,b)Q^{\mathsf T},\quad E(Qa,Qb)=QE(a,b)Q^{\mathsf T}.
\end{equation}
\end{proposition}
\begin{proof}
Expanding $M(a,b)b$ cancels the two terms proportional to $(a\cdot b)b$, leaving $\|b\|^2a$. The diagonal and transpose identities follow directly from \eqref{eq:M}. Equation \eqref{eq:cov} follows from $Q^{\mathsf T}Q=I$ and $(Qa)(Qb)^{\mathsf T}=Qab^{\mathsf T}Q^{\mathsf T}$.
\end{proof}

The symmetric and skew-symmetric parts of the product are respectively
\begin{equation}\label{eq:parts}
 \frac{M(a,b)+M(a,b)^{\mathsf T}}2=(a\cdot b)I,
\end{equation}
\begin{equation}
 \frac{M(a,b)-M(a,b)^{\mathsf T}}2=ab^{\mathsf T}-ba^{\mathsf T}.
\end{equation}
In general, $M(a,b)\ne M(b,a)$. The expression $M(M(a,b),c)$ lies outside the domain of the operation.

\section{Inverse vector and compatibility}
\begin{theorem}\label{thm:inverse}
For every $a\in\mathbb R^n\setminus\{0\}$, the equation $M(a,x)=I$ has exactly one solution $x\in\mathbb R^n$, namely
\begin{equation}\label{eq:inverse}
 a^{-1}=\frac{a}{\|a\|^2}.
\end{equation}
Moreover $M(a^{-1},a)=I$ and $(a^{-1})^{-1}=a$.
\end{theorem}
\begin{proof}
The notation $a^{-1}$ denotes only the unique vector solving $M(a,x)=I$; it does not make $\mathbb R^n$ a division algebra under $M$. Taking symmetric parts of $M(a,x)=I$ gives $a\cdot x=1$. Taking skew-symmetric parts gives $ax^{\mathsf T}-xa^{\mathsf T}=0$. Applying the latter matrix to $a$ yields $a(x\cdot a)-x\|a\|^2=0$, so $x=a/\|a\|^2$. Direct substitution verifies the equation and the reversed order. The double-inverse identity follows from $\|a^{-1}\|^2=\|a\|^{-2}$.
\end{proof}

\begin{corollary}\label{cor:compat}
For $b\ne0$,
\begin{equation}\label{eq:compat}
 E(a,b)=M(a,b^{-1}).
\end{equation}
\end{corollary}
\begin{proof}
Apply bilinearity of $M$ in its second argument and use \eqref{eq:inverse}.
\end{proof}

\section{An equivariant characterization}\label{sec:unique}
The transport property in isolation does not single out \eqref{eq:E}. Requiring a bilinear, orthogonally equivariant product with a prescribed scalar symmetric part does, however, give uniqueness.

\begin{theorem}\label{thm:unique}
Let $F:\mathbb R^n\times\mathbb R^n\to\mathbb R^{n\times n}$ be bilinear, $n\ge2$, and assume that, for all $a,b$ and $Q\in O(n)$,
\begin{equation}\label{eq:axioms}
 F(Qa,Qb)=QF(a,b)Q^{\mathsf T},
\end{equation}
\begin{equation}
 F(a,b)b=\|b\|^2a,\qquad \frac{F(a,b)+F(a,b)^{\mathsf T}}2=(a\cdot b)I.
\end{equation}
Then $F=M$.
\end{theorem}
\begin{proof}
By bilinearity, $F(a,b)_{ij}=T_{ijkl}a_kb_l$ for a fourth-order tensor $T$. Orthogonal equivariance implies that $T$ is invariant under $O(n)$. The standard classification of fourth-order orthogonally isotropic tensors \cite{mccullagh} yields
\begin{equation}\label{eq:form}
 F(a,b)=\alpha(a\cdot b)I+\beta ab^{\mathsf T}+\gamma ba^{\mathsf T}
\end{equation}
for constants $\alpha,\beta,\gamma$. For nonzero orthogonal $a,b$, the transport identity gives $\beta=1$. For arbitrary $a,b$ it then gives $\alpha+\gamma=0$. The prescribed symmetric part reads $\alpha(a\cdot b)I+\tfrac12(\beta+\gamma)(ab^{\mathsf T}+ba^{\mathsf T})=(a\cdot b)I$. Choosing nonzero orthogonal $a,b$ gives $\beta+\gamma=0$; choosing $a=b\ne0$ then gives $\alpha=1$. Thus $(\alpha,\beta,\gamma)=(1,1,-1)$, which is \eqref{eq:M}.
\end{proof}

\noindent\textbf{Role of the normalization hypotheses.} Transport alone does not ensure uniqueness: for each $t\in\mathbb R$, the bilinear, orthogonally equivariant map
\begin{equation}\label{eq:family}
 F_t(a,b)=t(a\cdot b)I+ab^{\mathsf T}-tba^{\mathsf T}
\end{equation}
satisfies $F_t(a,b)b=\|b\|^2a$. However, for $n\ge2$ its diagonal value is
\begin{equation}
 F_t(a,a)=t\|a\|^2I+(1-t)aa^{\mathsf T},
\end{equation}
which equals $\|a\|^2I$ for all $a$ only when $t=1$. Thus the diagonal identity is a genuine additional normalization, not a property shared by the family.

\begin{corollary}\label{cor:diagonal-unique}
Let $F$ be bilinear and $O(n)$-equivariant for $n\ge2$. If
$F(a,b)b=\|b\|^2a$ and $F(a,a)=\|a\|^2I$ for all $a,b$, then $F=M$.
\end{corollary}
\begin{proof}
Orthogonal equivariance gives the representation \eqref{eq:form}. Transport implies $\beta=1$ and $\alpha+\gamma=0$. On the diagonal,
$F(a,a)=\alpha\|a\|^2I+(1+\gamma)aa^{\mathsf T}$. For $n\ge2$, the matrices $I$ and $aa^{\mathsf T}$ are linearly independent for any $a\ne0$. The diagonal identity therefore implies $\alpha=1$ and $\gamma=-1$.
\end{proof}

The full orthogonal group is used in the isotropic-tensor argument; replacing it by a smaller symmetry group may permit additional equivariant tensors. Proposition~\ref{prop:without-equivariance} gives an explicit failure of uniqueness when equivariance is omitted.

\subsection{Why orthogonal equivariance cannot simply be dropped}
The diagonal normalization and transport identity alone do not
characterize the product in dimension three or higher. The following
explicit construction clarifies the role of the symmetry hypothesis.

\begin{proposition}\label{prop:without-equivariance}
Let $n\geq3$, fix a nonzero vector $c\in\mathbb R^n$, and let
$\omega:\mathbb R^n\times\mathbb R^n\times\mathbb R^n\to\mathbb R$
be a nonzero alternating trilinear form. For $\lambda\in\mathbb R$,
define the bilinear matrix-valued map $F_\lambda$ by
\begin{equation}\label{eq:counterexample-no-equiv}
 F_\lambda(a,b)x=M(a,b)x+\lambda\,\omega(a,b,x)c.
\end{equation}
Then $F_\lambda(a,b)b=\|b\|^2a$ and
$F_\lambda(a,a)=\|a\|^2I$ for all $a,b$.
For $\lambda\ne0$, $F_\lambda\ne M$.
\end{proposition}
\begin{proof}
For fixed $a,b$, the map $x\mapsto\omega(a,b,x)c$ is linear,
so it defines a matrix, and its dependence on $a,b$ is bilinear.
Alternation gives $\omega(a,b,b)=0$ and
$\omega(a,a,x)=0$. The two stated identities therefore follow
from the corresponding identities for $M$. Since $\omega$ is nonzero
and $c\ne0$, there exist $a,b,x$ with $\omega(a,b,x)c\ne0$;
thus $F_\lambda\ne M$ whenever $\lambda\ne0$.
\end{proof}

For example, in $\mathbb R^3$ take
$\omega(a,b,x)=\det[a\ b\ x]$ and $c=e_1$.
The additional term in \eqref{eq:counterexample-no-equiv} is then
$\lambda e_1(a\times b)^{\mathsf T}x$.
It selects a fixed spatial direction $e_1$ and hence does not obey
full orthogonal equivariance. This counterexample shows that
bilinearity, transport, and diagonal normalization do not suffice
without an appropriate symmetry condition for $n\ge3$.

\subsection{The special case of dimension two}
The preceding counterexample requires dimension at least three. In the
plane, transport and diagonal normalization already suffice without
assuming orthogonal equivariance.

\begin{proposition}\label{prop:two-dimensional-unique}
Let $F:\mathbb R^2\times\mathbb R^2\to\mathbb R^{2\times2}$ be
bilinear. If $F(a,b)b=\|b\|^2a$ and
$F(a,a)=\|a\|^2I_2$ for every $a,b$, then $F=M$.
\end{proposition}
\begin{proof}
Set $D=F-M$. Then $D$ is bilinear, $D(a,a)=0$, and
$D(a,b)b=0$. Polarizing $D(a,a)=0$ gives
$D(a,b)=-D(b,a)$. With the standard basis $e_1,e_2$, every
alternating bilinear matrix-valued map on $\mathbb R^2$ has the form
$D(a,b)=\det[a\ b]C$ for a fixed $2\times2$ matrix $C$.
The transport identity implies $Cb=0$ whenever $b\ne0$ (choose
$a$ independent of $b$). Thus $C=0$ and $F=M$.
\end{proof}

\section{A constrained minimum-norm characterization}\label{sec:min}
Set
\begin{equation}\label{eq:alpha}
 \alpha=\frac{a\cdot b}{\|b\|^2},\qquad
 K_0=\frac{ab^{\mathsf T}-ba^{\mathsf T}}{\|b\|^2}.
\end{equation}
Then $E(a,b)=\alpha I+K_0$, with $K_0^{\mathsf T}=-K_0$.

\begin{theorem}\label{thm:min}
Among matrices $X\in\mathbb R^{n\times n}$ satisfying $Xb=a$ and having scalar symmetric part, that is,
\begin{equation}\label{eq:sym}
 \frac{X+X^{\mathsf T}}2=sI\quad\text{for some }s\in\mathbb R,
\end{equation}
the unique minimizer of the Frobenius norm is $E(a,b)$.
\end{theorem}
\begin{proof}
The scalar is forced to be $s=\alpha$ by $b^{\mathsf T}Xb=a\cdot b$. Any feasible $X$ has the form $E(a,b)+H$ with $H^{\mathsf T}=-H$ and $Hb=0$. Since $b^{\mathsf T}H=0$ as well,
\begin{equation}
 \langle K_0,H\rangle_F
 =\frac{\operatorname{tr}[(ba^{\mathsf T}-ab^{\mathsf T})H]}{\|b\|^2}=0.
\end{equation}
The scalar matrix is Frobenius-orthogonal to every skew-symmetric matrix. Consequently,
\begin{equation}
 \|X\|_F^2=\|E(a,b)\|_F^2+\|H\|_F^2,
\end{equation}
with equality at the minimum only for $H=0$.
\end{proof}

Direct computation gives
\begin{equation}\label{eq:frob}
 \|E(a,b)\|_F^2=2\frac{\|a\|^2}{\|b\|^2}+(n-2)\alpha^2.
\end{equation}
The symmetric-part constraint is indispensable: without it, the minimum-Frobenius-norm solution of $Xb=a$ is the rank-one matrix $ab^{\mathsf T}/\|b\|^2$.

\section{Spectral structure and limits of reciprocity}\label{sec:spectrum}
Assume first that $a,b$ are nonzero and linearly independent. Let $P=\operatorname{span}\{a,b\}$,
\begin{equation}
 \rho=\frac{\|a\|}{\|b\|},\qquad
 \theta=\arccos\!\left(\frac{a\cdot b}{\|a\|\|b\|}\right)\in(0,\pi).
\end{equation}
Writing $u=a/\|a\|$, $v=b/\|b\|$ and choosing $e_1=v$, $e_2=(u-\cos\theta\,v)/\sin\theta$, the restriction of $E(a,b)$ to $P$ has matrix $\rho\left(\begin{smallmatrix}\cos\theta&-\sin\theta\\\sin\theta&\cos\theta\end{smallmatrix}\right)=\rho R_\theta$. On $P^\perp$ it acts as the scalar $\rho\cos\theta$. Therefore,
\begin{equation}\label{eq:decomp}
 E(a,b)\simeq\rho R_\theta\oplus\rho\cos\theta\,I_{n-2}.
\end{equation}
Its complex eigenvalues are $\rho e^{\pm i\theta}$ and $\rho\cos\theta$ with multiplicity $n-2$. In particular,
\begin{equation}\label{eq:det}
 \det E(a,b)=\begin{cases}\dfrac{\|a\|^2}{\|b\|^2},&n=2,\\[3pt]\dfrac{\|a\|^2}{\|b\|^2}\left(\dfrac{a\cdot b}{\|b\|^2}\right)^{n-2},&n>2.\end{cases}
\end{equation}
The determinant formula extends to nonzero collinear pairs by continuity, where $E(a,b)$ is a scalar multiple of $I$. If $a=0$, then $E(0,b)=0$. For $a\ne0$, the quotient matrix is invertible in dimension two; for $n>2$, it is invertible precisely when $a\cdot b\ne0$.

\begin{proposition}\label{prop:recip}
For nonzero linearly independent $a,b$, let $\Pi_{P^\perp}$ be orthogonal projection onto $P^\perp$. Then
\begin{equation}\label{eq:recip}
 E(a,b)E(b,a)=I-\sin^2\theta\,\Pi_{P^\perp}.
\end{equation}
\end{proposition}
\begin{proof}
The factors restrict on $P$ to $\rho R_\theta$ and $\rho^{-1}R_{-\theta}$, whose product is $I_P$. On $P^\perp$ their scalar factors multiply to $\cos^2\theta$.
\end{proof}

For nonzero collinear pairs the product is $I$. Thus $E(a,b)E(b,a)=I$ for all nonzero $a,b$ in dimension two, whereas for $n>2$ it holds if and only if $a,b$ are collinear. The vector inverse in Theorem~\ref{thm:inverse} should not be confused with the inverse of the matrix $E(a,b)$.

\section{Explicit examples}
\paragraph{Dimension two.}
For $a=(3,-1)^{\mathsf T}$ and $b=(2,5)^{\mathsf T}$, $a\cdot b=1$ and $\|b\|^2=29$. Hence
\begin{equation}
 M(a,b)=\begin{pmatrix}1&17\\-17&1\end{pmatrix},\qquad
 E(a,b)=\frac1{29}\begin{pmatrix}1&17\\-17&1\end{pmatrix}.
\end{equation}
Direct multiplication gives $E(a,b)b=a$, while $a^{-1}=a/10$ satisfies $M(a,a^{-1})=I_2$.

\paragraph{Dimension three.}
For $a=(3,-1,2)^{\mathsf T}$ and $b=(2,5,1)^{\mathsf T}$, $a\cdot b=3$ and $\|b\|^2=30$, so
\begin{equation}
 M(a,b)=\begin{pmatrix}3&17&-1\\-17&3&-11\\1&11&3\end{pmatrix},
 \qquad E(a,b)=\frac{M(a,b)}{30}.
\end{equation}
Thus $E(a,b)b=a$ and $a^{-1}=a/14$ gives $M(a,a^{-1})=I_3$.

\paragraph{Singularity despite transport.}
For the standard basis vectors $e_1,e_2\in\mathbb R^3$,
\begin{equation}
 E(e_1,e_2)=\begin{pmatrix}0&1&0\\-1&0&0\\0&0&0\end{pmatrix}.
\end{equation}
This matrix maps $e_2$ to $e_1$ but has determinant zero.

\section{Relation to established algebraic structures}\label{sec:context}
The decomposition in \eqref{eq:parts} separates the scalar inner product from the skew operator $ab^{\mathsf T}-ba^{\mathsf T}$. The latter is naturally associated with the bivector $a\wedge b$. In Clifford (geometric) algebra, the geometric product likewise contains scalar and bivector contributions; however, the present product has matrix codomain and is not identified with the associative Clifford product \cite{hestenes,clifford}. In particular, $a/\|a\|^2$ is also the standard inverse of a non-null Euclidean vector in geometric algebra, so the formula for the inverse alone is not a novelty claim.

In dimension two, with $z_a=a_1+ia_2$, $z_b=b_1+ib_2$ and $\det[a\ b]=a_1b_2-a_2b_1$, the complex quotient is $z_a/z_b=((a\cdot b)-i\det[a\ b])/\|b\|^2$; its real multiplication matrix is precisely $E(a,b)$. Under this convention, $M(a,b)$ corresponds to $z_a\overline{z_b}$, rather than $z_az_b$. The two-dimensional quotient is therefore an established complex-algebra operation in matrix form. In dimension three, with $[w]_\times x=w\times x$, one has $ab^{\mathsf T}-ba^{\mathsf T}=-[a\times b]_\times$. Equations \eqref{eq:E}--\eqref{eq:M} nevertheless provide one dimension-independent matrix formulation, together with the specific equivariant and variational characterizations established above.

\subsection{Comparison with rotor-induced rotations}\label{sec:rotor}
Let $a,b$ be nonzero and linearly independent, and put
$u=a/\|a\|$, $v=b/\|b\|$, $P=\operatorname{span}\{a,b\}$,
$\rho=\|a\|/\|b\|$, and $c=u\cdot v$.
Denote by $\Pi_{P^\perp}$ the orthogonal projector onto $P^\perp$.
The orientation-preserving orthogonal map $R_{u\leftarrow v}$ that
maps $v$ to $u$ and fixes $P^\perp$ pointwise is the minimal planar
rotation between the two directions.
\begin{proposition}\label{prop:rotor}
With the notation above,
\begin{equation}\label{eq:rotor-relation}
 E(a,b)=\rho R_{u\leftarrow v}-\rho(1-c)\Pi_{P^\perp}.
\end{equation}
In particular, $E(a,b)|_P=\rho R_{u\leftarrow v}|_P$ and
$E(a,b)|_{P^\perp}=\rho c I_{P^\perp}$.
\end{proposition}
\begin{proof}
On $P$, the sum of the scalar and skew-symmetric terms in $E(a,b)$ equals $\rho R_{u\leftarrow v}|_P$ (as is seen in the adapted orthonormal basis of Section~\ref{sec:spectrum}). On $P^\perp$, both rank-one terms vanish, and the
remaining scalar is $(a\cdot b)/\|b\|^2=\rho c$.
The right-hand side of \eqref{eq:rotor-relation} has precisely
these restrictions, proving equality.
\end{proof}
For the Euclidean Clifford product, the normalized even element
\begin{equation}\label{eq:rotor}
 U=\frac{1+uv}{\sqrt{2(1+u\cdot v)}}
\end{equation}
is a rotor, since $u\ne-v$, and its sandwich action satisfies
$U v\widetilde U=u$, where $\widetilde U$ denotes reversion.
The corresponding vector-space action is
$R_{u\leftarrow v}x=Ux\widetilde U$; thus
\begin{equation}\label{eq:clifford-operator}
 E(a,b)x=\rho Ux\widetilde U
 -\rho(1-c)\Pi_{P^\perp}x.
\end{equation}
Here $U$ is a Clifford-algebra element, whereas $E$ and $R$ are
real-linear operators. Equation \eqref{eq:clifford-operator} does not
identify the matrix-valued product $M$ with Clifford multiplication.
Indeed, the linear identification of a bivector $a\wedge b$ with
the skew endomorphism $ab^{\mathsf T}-ba^{\mathsf T}$ is not a
multiplicative representation of the full Clifford algebra.
For example, for orthonormal $e_1,e_2$ in $\mathbb R^3$,
$(e_1e_2)^2=-1$ in Euclidean Clifford algebra, while
$M(e_1,e_2)^2=-\Pi_{\operatorname{span}\{e_1,e_2\}}\ne-I_3$.
For $u=-v$ and $n\ge3$, a planar half-turn requires an additional choice of a two-dimensional plane containing $v$. In dimension two, the corresponding rotation is uniquely $-I_2$. See Hestenes' treatment of vectors, spinors, and
rotations \cite{hestenes1971,hestenes}.

\section{Conclusions}
The matrix-valued quotient \eqref{eq:E} satisfies an exact transport identity and is compatible with the bilinear product \eqref{eq:M}. The identity equation $M(a,x)=I$ determines the unique inverse vector $x=a/\|a\|^2$, yielding $E(a,b)=M(a,b^{-1})$. Orthogonal equivariance and transport, together with either the prescribed symmetric part or the diagonal identity, characterize the bilinear product; an explicit alternating-tensor counterexample shows why the symmetry hypothesis matters in dimensions at least three, whereas dimension two needs no such assumption; a constrained Frobenius-norm principle selects the quotient. The spectral decomposition specifies its invertibility and explains the dimension-dependent failure of matrix reciprocity. An explicit correction on the orthogonal complement relates the quotient to a rotor-induced planar rotation, without identifying the matrix product with Clifford multiplication.

The contribution is the self-contained synthesis of these matrix-valued operations and their axiomatic, variational, and spectral properties. No priority claim is made for the underlying classical algebraic constructions or for the general theory of isotropic tensors. Extensions to complex vectors and composition along multi-step paths are beyond the scope of this paper.

\end{document}